\documentclass[12 pt]{article}
\usepackage{amsfonts}

\usepackage{graphicx}
\usepackage{amsmath}
\usepackage{amssymb}

\catcode`\Š=\active\defŠ{\"a}        
\catcode`\'=\active\def'{\"e}        
\catcode`\•=\active\def•{\"{\i}}     
\catcode`\š=\active\defš{\"o}        
\catcode`\Ÿ=\active\defŸ{\"u}        
\catcode`\Ø=\active\defØ{\"y}        
\catcode`\€=\active\def€{\"A}        
\catcode`\…=\active\def…{\"O}        
\catcode`\†=\active\def†{\"U}        
\catcode`\‡=\active\def‡{\'a}        
\catcode`\Ž=\active\defŽ{\'e}        
\catcode`\'=\active\def'{\'{\i}}     
\catcode`\—=\active\def—{\'o}        
\catcode`\œ=\active\defœ{\'u}        
\catcode`\ƒ=\active\defƒ{\'E}        
\catcode`\ˆ=\active\defˆ{\`a}        
\catcode`\=\active\def{\`e}        
\catcode`\"=\active\def"{\`{\i}}     
\catcode`\˜=\active\def˜{\`o}        
\catcode`\=\active\def{\`u}        
\catcode`\Ë=\active\defË{\`A}        
\catcode`\‹=\active\def‹{\~a}        
\catcode`\–=\active\def–{\~n}        
\catcode`\›=\active\def›{\~o}        
\catcode`\Ì=\active\defÌ{\~A}        
\catcode`\"=\active\def"{\~N}        
\catcode`\Í=\active\defÍ{\~O}        
\catcode`\‰=\active\def‰{\^a}        
\catcode`\=\active\def{\^e}        
\catcode`\"=\active\def"{\^{\i}}     
\catcode`\™=\active\def™{\^o}        
\catcode`\ž=\active\defž{\^u}        

\begin{document}
\large
\begin{center}
\textbf{Expansions of the Riemann Zeta function  in the  critical strip}
\medskip

 \text{ B. Candelpergher }\footnote{ Univ. de Nice, Lab. J.A.DieudonnŽ , email: candel@unice.fr}

\end{center}

\bigskip

\textbf{Introduction}

\bigskip

We introduce the functions defined for $ t\in ]0,+\infty[
$ by 

$${
\Psi_m(t)={\sqrt{2}}\ (\frac{t^2-1}{t^2+1})^m\frac1{\sqrt{1+t^2}}}
$$

where $m\geq 0$ is an integer. Their  Mellin transform are 
 
 $$
\mathcal{M}(\Psi_m)(s)= \int_0^{+\infty} t^{ s-1} {\Psi_m(t)}dt
= \frac1{\sqrt{2\pi}} \Gamma(\frac{s}2)\Gamma(\frac{1-s}2)Q_{m}(s)$$

 where  $Q_m$ are polynomials in $\mathbb{R}[X] $ with  their roots on the line  $Re(s)=1/2$.
 
We use these functions $\Psi_m$  to get the expansion  

$$
\sum_{n\in\mathbb{Z}}e^{-\pi n^2 t^2}-1-\frac1t= \sum_{m\geq 0}{\alpha_{2m}}\Psi_{2m}(t)
Ê\ \text{ for }t\in ]0,+\infty[
$$

\text{ with }

$$
\alpha_{2m}=\frac{2^{-4m}}{(2m)!}\Big(\sum_{n\in \mathbb{Z}}
H_{4m}(\sqrt{2\pi}\ n)e^{-\pi n^2}-2\frac{(4m)!}{(2m)!}\Big)
$$

where $H_n$ are the Hermite polynomials.

In the strip $0<Re(s)<1$ the well-known classical result 
$$
\int_0^{+\infty} t^{s-1}(\sum_{n\in\mathbb{Z}}e^{-\pi n^2 t^2}-1-\frac1t)dt
=\Gamma(\frac s2)\pi^{-s/2}\zeta(s)
$$

allows us to conjecture the following expansion of Zeta for $0<Re(s)<1$
$${
\zeta(s)=\frac1{\sqrt{2\pi}}\ \pi^{\frac{s}2}\ \Gamma(\frac{1-s}2) \sum_{m\geq 0}{\alpha_{2m}}\ Q_{2m}(s)}
$$

\section {Functions related to the quantum harmonic oscillator}

\subsection{Hermite functions}

For every integer  $m\geq 0$ let us consider the Hermite function
 $$
\Phi_m (x)=H_m(\sqrt{2\pi} x)e^{-\pi x^2}
$$

where $H_m\in\mathbb{R}[X] $ are the Hermite polynomials defined  by the generating function 
$$
e^{-t^2+2xt}=\sum_{m\geq 0}\frac{H_m(x)}{m!}t^m
$$

or directly by 
$H_m(x)=(-1)^me^{x^2}\partial^m e^{-x^2}.$
 
The Hermite functions $\Phi_m\in L^2(\mathbb{R})$ are known to form an orthogonal system of eigenfunctions of the quantum harmonic oscillator
$$
2\pi(x^2-\frac1{4\pi^2}\partial^2)\Phi_m=(2m+1)\Phi_m 
\ \text{ with }\ \int_{\mathbb{R}}(\Phi_m(x))^2dx=\frac1{\sqrt{2}}2^mm!
$$

 The function $\Phi_m $ has same parity as $m$.
We are  interested with the even functions $\Phi_{2m}$,  we have (cf. [4]) 
$$
\int_{\mathbb{R}}e^{-2i\pi x\xi}\Phi_{2m}(x)dx=(-1)^m\Phi_{2m}(\xi)
$$

Thus for $\xi=0$
$$
\int_{\mathbb{R}}\Phi_{2m}(x)dx=(-1)^m\Phi_{2m}(0)=\frac{(2m)!}{m!}
$$

The function $\Phi_{2m}$ is bounded  (cf. [4])  by
$$
B1)\hspace{1cm}Ê
\vert \Phi_{2m}(x)\vert\leq K2^{m} \sqrt{(2m)!} \ \text{ with }K=1.086435 \ \text{ for } x \in \mathbb{R}\\
$$

The function $\Phi_{2m}$ is (cf. [5]) oscillating in the interval 
$$I_m=[-\frac1{\sqrt{\pi}}\sqrt{2m+1}, \frac1{\sqrt{\pi}}\sqrt{2m+1}]$$

 and  exponentially decreasing when $x\notin I_m$
 , more precisely (cf. [4]) we have 
\begin{eqnarray*}
(B2)\hspace{1cm}Ê\vert \Phi_{2m}(x)\vert\leq \frac{(2m)!}{m!} e^{2x\sqrt{2\pi m}}e^{-\pi x^2}  \text{ for } x>0
\end{eqnarray*}

\medskip

In the  series expansions of the following sections we use the normalized sums $$S_{2m}=\frac{2^{-2m}}{m!}\sum_{n\in \mathbb{Z}} \Phi_{2m}(n)$$

\bigskip

\textbf{Lemma 0}

For $m\rightarrow +\infty$ we have
$$
\frac{2^{-2m}}{m!}\sum_{ \vert n\vert \geq 2\sqrt{2m} } \vert\Phi_{2m}(n)\vert =  O(e^{-2\pi\sqrt{2m}})
$$

and
$$\frac{2^{-2m}}{m!}\sum_{n\in\mathbb{Z}}\vert \Phi_{2m}(n)\vert=O(m^{1/4})$$

\bigskip

\textit{Proof}

We have for  $m\geq 1$ 
$$
2x\sqrt{2\pi m}-\pi x^2\leq  -\pi   x \ \text{ for }\  x\geq 2\sqrt{2m}
$$

 thus using  inequality (B2)  we get  
$$
\vert \frac{2^{-2m}}{m!}\Phi_{2m}(x)\vert\leq \frac{2^{-2m}(2m)!}{(m!)^2} e^{-\pi  \vert x \vert} \text{ for }
\vert x \vert \geq 2\sqrt{2m}$$

Thus by summation for $\vert n \vert \geq 2\sqrt{2m}$ and with the Stirling formula we get 
 \begin{eqnarray*}
 \frac{2^{-2m}}{m!}\sum_{\vert n\vert\geq 2\sqrt{2m}}\vert \Phi_{2m}(n)\vert\leq
\frac{2^{-2m+1}(2m)!}{(m!)^2}\frac{e^{-\pi (2\sqrt{2m})}
}{1-e^{-\pi }} =O(e^{-\pi 2\sqrt{2m}})\end{eqnarray*}

 From inequality (B1) we deduce that
  \begin{eqnarray*}
 \frac{2^{-2m}}{m!}\sum_{n\in\mathbb{Z}}\vert \Phi_{2m}(n)\vert&=&\sum_{\vert n\vert < 2\sqrt{2m}}\frac{2^{-2m}}{m!}\vert \Phi_{2m}(n)\vert+\frac{2^{-2m}}{m!}\sum_{\vert n\vert\geq 2\sqrt{2m}}\vert \Phi_{2m}(n)\vert\\
 &\leq&  K\sqrt{2m}\frac{2^{-m+1}\sqrt{(2m)!}}{m!}+O(e^{-\pi 2\sqrt{2m}}) \end{eqnarray*}
 
 thus by Stirling formula we get 
  $\frac{2^{-2m}}{m!}\sum_{n\in\mathbb{Z}}\vert \Phi_{2m}(n)\vert=O(m^{1/4}).$
  
  $\square$

\subsection{The functions $\Psi_m$}
The function $x\mapsto e^{-2\pi a^2x^2}$, $Re(a^2)>-\frac12$, expands (cf. [6] p.71-75) 
 as the following series of Hermite polynomials,   for $x\in \mathbb{R}$ we have

$$
e^{-2\pi a^2x^2}=\frac1{\sqrt{1+a^2}}\sum_{m\geq 0}\frac{(-1)^ma^{2m}}{2^{2m}(1+a^2)^mm!}H_{2m}(\sqrt{2\pi}\ x)
$$

Multiplying by $e^{-\pi x^2}$ we get, with    $t^2=1+2a^2$
\begin{equation}
e^{-\pi x^2 t^2}=\sum_{m\geq 0}{(-1)^m}\frac1{2^{2m}}\Phi_{2m}(x)\frac{\Psi_m(t)}{m!} \text{  for }Re(t^2)>0
\end{equation}

where we define for $t\in S=\{re^{i\theta}\vert\  r>0,-\frac{\pi}4<\theta<\frac{\pi}4\}$ the function
$${
\Psi_m(t)={\sqrt{2}}\ (\frac{t^2-1}{t^2+1})^m\frac1{\sqrt{1+t^2}}}
$$

\medskip

\textbf{Lemma 1}

 The functions $\Psi_m$ are related to the Hermite functions by
$${
\frac {\Psi_m}{m!}=(\frac{2\sqrt{2}}xe^{-\pi /x^2})\ * 
 \frac{\Phi_{2m}}{(2m)!} }
$$

where $*$ is  the multiplicative convolution  of functions defined on $]0,+\infty[$ 
 $$
(f*g)(t)=\int_0^{+\infty}f(\frac tx)g(x)\frac1x dx$$

 \bigskip
 
 \textit{Proof}
 
 With the classical relation
$$
 \int_{-\infty}
^{+\infty}
 e^{-a x^2}e^{b x}dx=\sqrt{\frac{\pi }a}\ e^{\frac{b^2}{4a}}\ \text{ where } a>0,\ b\in \mathbb{C}
 $$
 
we get 
$$
e^{z^2\frac{t^2-1}{t^2+1}}\frac t{\sqrt{1+t^2}}=\int_{-\infty}
^{+\infty}
 e^{-\pi x^2\frac {t^2+1}{t^2}}e^{-z^2+2\sqrt{2\pi} xz}dx
 $$
 
and using the power series expansion
 $$
 e^{-z^2+2\sqrt{2\pi} xz}=\sum_{m\geq 0}\frac{z^m}{m!}H_{m}(\sqrt{2\pi} x) 
 $$
 
 we get by identification
 $$
\frac 1{m!}\sqrt{2}\ (\frac{t^2-1}{t^2+1})^m\frac t{\sqrt{1+t^2}}=2\sqrt{2}\int_0^{+\infty}
 e^{-\pi x^2\frac {t^2+1}{t^2}}
\frac{1}{(2m)!}H_{2m}(\sqrt{2\pi} x) dx
$$

This gives
 $${
\frac 1{m!}\Psi_m(t)=2\sqrt{2}\int_0^{+\infty}
 e^{-\pi x^2/t^2}\frac1t
\frac{\Phi_{2m}(x)}{(2m)!} dx}
$$
 
and we see that this last  integral  is the multiplicative convolution  
 $$
(f*g)(t)=\int_0^{+\infty}f(\frac tx)g(x)\frac1x dx
$$

 with  $f(x)=2\sqrt{2}\ e^{-\pi /x^2}\frac1x$ 
  and $g(x)=\frac{\Phi_{2m}(x)}{(2m)!} $.

$\square$

\section{Series expansions}

Let  $z\mapsto \sqrt{z}$ the principal determination of the square root, the holomorphic function
$$
u\mapsto t=\sqrt{\frac{1+u}{1-u} }
$$

maps the open unit disk $D(0,1)=\{z\in \mathbb{C}\vert\  \ \vert z\vert <1\}$ onto the sector
$$S=\{re^{i\theta}\vert\  r>0,-\frac{\pi}4<\theta<\frac{\pi}4\}$$ 

For any function $f$ holomorphic in  $S$ let us define the function

$${
Tf(u)=\frac{1}{\sqrt{1-u}}f\Big(\sqrt{\frac{1+u}{1-u}}\Big)}
$$

which is holomorphic in the open disk  $D(0,1)$.

For every integer $m\geq 0$  we verify immediately that we have
 $$T\Psi_m(u)=u^m$$

For a function  $f$ defined on  $S$ the expansion
$$
f(t)=\sum_{m\geq 0}a_m\frac {\Psi_m(t)}{m!}$$

follows the Taylor expansion of  $Tf$ 
$$
Tf(u)=\sum_{m\geq 0}\frac {a_m}{m!}\ u^m
$$

\medskip

\textbf{Remark.} 
Note that $\Psi_m(t)=(-1)^m\frac1t\Psi_m(\frac1t)$ for all $t\in S$. 
For a function $f$ on  $S$   the relation   $$f(t)=\frac1tf(\frac1t)$$ 

is equivalent to the parity of the function $Tf$
$$Tf(u)=Tf(-u)$$

in this case the expansion of $f$ is of the form
$$
f(t)=\sum_{m\geq 0}a_{2m}\frac {\Psi_{2m}(t)}{(2m)!}$$

\medskip

\textbf{Example}

For the function $f=\frac1{1+t}$, $t\in S$, we have 
 $$
Tf(u)=\frac{\sqrt{1+u}-\sqrt{1-u}}{2u}
$$

This gives for $t\in S$
 
\begin{equation}
 \frac1{1+t}=\frac1{{2}}\sum_{m\geq 0} \frac{(4m)!}{2^{4m}(2m+1)!}\frac {\Psi_{2m}(t)}{(2m)!}
 \end{equation}
 
\bigskip

\subsection{Expansion of the theta function}

The theta function defined for $t\in S$ by
$$
G(t)=\sum_{n\in\mathbb{Z}}e^{-\pi n^2 t^2}$$
 
is holomorphic in $S$  and we have for  $u\in D(0,1)$

$$TG(u)=\frac{1}{\sqrt{1-u}}\sum_{n\in\mathbb{Z}}e^{-\pi n^2 \frac{1+u}{1-u}}$$

Let 
$$TG(u)=\sum_{m\geq 0 } g_m \frac1{m!}u^m$$

be the power series expansion of the holomorphic function $TG$ in the open disk $D(0,1)$.
The  Jacobi identity  (cf.[3]) 
$$\frac1tG(\frac1t)= G(t)$$

gives the parity of  $TG$ and we get
$TG(u)=\sum_{n\geq 0 } g_{2m} \frac1{(2m)!}u^{2m}$.

Thus we have  for $t\in S$
 
 $$G(t)=\sum_{m\geq 0}g_{2m}\frac{\Psi_{2m}(t)}{(2m)!}$$
 
 \bigskip

\textbf{Lemma 2}

We have for $t\in S=\{re^{i\theta}\vert \  r>0,-\frac{\pi}4<\theta<\frac{\pi}4\}$

$$
G(t)=\sum_{m\geq 0}S_{4m}{\Psi_{2m}(t)}
\  \text{ where } \ S_{4m}=\frac{2^{-4m}}{(2m)!}
\sum_{n\in \mathbb{Z}} \Phi_{4m}(n)$$

\bigskip

\textit{Proof}

Take the relation (1) with $x=n\in\mathbb{Z}$, by summation we get 

\begin{eqnarray*}
\sum_{n\in\mathbb{Z}}e^{-\pi n^2 t^2}=\sum_{n\in\mathbb{Z}}\sum_{m\geq 0}\frac{(-1)^m}{m!}2^{-2m}\Phi_{2m}(n)\Psi_m(t)\\
=\sum_{m\geq 0}\frac{(-1)^m}{m!}2^{-2m}\Psi_m(t)\sum_{n\in\mathbb{Z}}\Phi_{2m}(n)
\end{eqnarray*}

To justify the interchange of summations $\sum_{n\in\mathbb{Z}}\sum_{m\geq 0}=\sum_{m\geq 0}\sum_{n\in\mathbb{Z}}$ we observe that
 $$\vert \frac{t^2-1}{t^2+1}\vert<1 \text { for   }t\in S
 $$
 and by Lemma 0 we have $ \frac{2^{-2m}}{m!}\sum_{n\in\mathbb{Z}}\vert \Phi_{2m}(n)\vert=O(m^{1/4})$ 
 thus  for $t\in S$
 
 $$
\sum_{m\geq 0}\frac{2^{-2m}}{m!}\sum_{n\in\mathbb{Z}}\vert \Phi_{2m}(n)\vert \vert \Psi_m(t)\vert<+\infty
$$

This gives 
$$
G(t)=\sum_{m\geq 0}{(-1)^m}S_{2m}\Psi_m(t)
$$

Since we have seen that the function $TG$ is even
we deduce that in this last sum, only the constants  $S_{4m}$ are non  zero. 

$\square$

\textbf{Remark}

We have  also (cf. Appendix)  for  the constants $S_{4m}$ another expression
$$
S_{4m}=
\frac{(\frac{\pi}2)^{2m}}{(2m)!}
\frac1{S_0\sqrt{2}}\sum_{(k,l)\in  \mathbb{Z}^2}(-1)^{kl}e^{-\pi\frac12(k^2+l^2)}(k+il)^{4m}
$$

\medskip

\textbf{Theorem}

For $t\in S=\{re^{i\theta}\vert \  r>0,-\frac{\pi}4<\theta<\frac{\pi}4\}$ we have
$${
G(t)-1-\frac1t= \sum_{m\geq 0}{\alpha_{2m}}\Psi_{2m}(t)
\text{ with }
\alpha_{2m}=S_{4m}-\frac{{2^{-4m+1}}(4m)!}{(2m)!(2m)!}}
$$

\textit{Proof}

Using Lemma 2, to get the expansion of $G(t)-1-\frac1{t}$ in terms of $\Psi_m(t)$ it is now sufficient to  expand  $1+\frac1t$.
For $f(t)=1+\frac1t$ one has 
$$
Tf(u)=\frac{1}{\sqrt{1+u}}+\frac{1}{\sqrt{1-u}}
$$

and we obtain for $t\in S$

$$
1+\frac1t= {2}\sum_{m\geq 0}\frac{(4m)!}{2^{4m}(2m)!}\frac{\Psi_{2m}(t)}{(2m)!}
$$
$\square$

\pagebreak

\textbf{Remark}

We see that 
$$\alpha_{2m}=\frac{2^{-4m}}{(2m)!}\Big(\sum_{n\neq 0}\Phi_{4m}(n)-[\Phi_{4m}(0)+\int_{\mathbb{R}}\Phi_{4m}(x)dx]\Big)$$

This is easily explained if we look at the general MŸntz formula (cf. [7]):

 let $F$ be an even continuously differentiable function such that $F$ and $F'$  are $O(x^{-a})$, $(a>1)$ when $x\rightarrow \infty$, then for $0<Re(s)<1$ we have
$$
2\ \zeta(s)\mathcal{M}F(s)=\mathcal{M}\Big(G(t)-[F(0)+\int_{\mathbb{R}}F(xt)dx]\Big)(s)
$$

with $G(t)=\sum_{n\in \mathbb{Z}} F(nt)$, this is our case with $F(x)=e^{-\pi x^2}$.

If there exist functions a sequence of functions $\varphi_m$ and $\psi_m$ such that we have an expansion
$$
F(xt)=\sum_{m\geq 0} \varphi_m(x)\psi_m(t)
$$

then, at least formally, we get 
$$
G(t)-[F(0)+\int_{\mathbb{R}}F(xt)dx]=\sum_{m\geq 0}\Big(\sum_{n\in \mathbb{Z}}\varphi_m(n)-[\varphi_m(0)+\int_{\mathbb{R}}\varphi_m(x)dx]\Big)\psi_m(t)
$$

in our case $\varphi_m=\frac{2^{-4m}}{(2m)!}\Phi_{4m}$ and 
$\psi_m(t)=\Psi_{2m}(t)$.



\section{ Mellin transforms}

\subsection{The   polynomials $Q_m$}

For $Re(s)>0$,  the Mellin transforms of the Hermite  functions $\Phi_{2m}$ are

$${\int_0^{+\infty} \frac{\Phi_{2m}(x)}{(2m)!}x^{s-1}dx=\frac12\pi^{-s/2}\Gamma(s/2)\frac{Q_m(s)}{m!}}
$$

where $Q_m$  are polynomials in $\mathbb{R}[X] $. 
This is simply a consequence of the relation  
$$
\int_0^{+\infty} e^{-\pi x^2}x^{s+2k-1}dx
=\frac12\pi^{-\frac s2}\Gamma(\frac s2)
\ \pi^{-k}\ \frac s2(\frac s2+1)...(\frac s2+k-1)$$

We get $Q_0(s)=1$, 
$Q_1(s)={2s-1}$,
$Q_2(s)={\frac43s^2-\frac43 s+1}$,... .

\bigskip

 More generally we have
 $$
Q_{m}(s)=\sum_{k=0}^m(-1)^{m-k}\frac{m!}{(m-k)!}\frac{2^{2k}}{(2k)!}s(s+2)...(s+2(k-1))
$$

and (cf. [4]) an expression of $Q_m$ in terms of the hypergeometric function
$$
Q_{m}(s)={(-1)^m}\ _2F_1(-m,s/2;1/2;2)
$$

\medskip

The roots of  $Q_m$ are  on the line $Re(s)=1/2$ (cf. [1], [2]). This can be proved (cf. [1]) by observing that  the orthogonality relation of the Hermite functions $\Phi_{2m}$ implies the orthogonality of the family of polynomials $$t\mapsto Q_m(\frac12+it)$$ with respect to the Borel measure $\vert \Gamma(\frac14+i\frac t2)\vert^2 dt $ on $\mathbb{R}$.

 More explicitly, using the Parseval's formula for Mellin transform
 $$
 \frac1{2\pi}\int_{-\infty}^{+\infty}\big(\mathcal{M}(f)\overline{\mathcal{M}(g)}\big)(\frac12+it)dt=\int_{0}^{+\infty}(f\overline{g})(x)dx
$$

  we get 
  $$
\frac1{4\pi\sqrt{\pi}}\int_{\mathbb{R}}\vert \Gamma(\frac14+i\frac t2)\vert^2\ \big(\frac{Q_{m_1}}{m_1!}
\frac{\overline{Q_{m_2}}}{m_2!}\big)(\frac12+it)dt=\int_{\mathbb{R}} \big(\frac{\Phi_{2m_1}}{(2m_1)!}\frac{\Phi_{2m_2}}{(2m_2)!}\big)(x)dx
$$

\subsection {Mellin transform of $\Psi_m$}

\textbf{Lemma 3}

For $0<Re(s)<1$ we have
$$
{\int_0^{+\infty} t^{ s-1}\Psi_m(t)dt
= \frac1{\sqrt{2\pi}} \Gamma(\frac{s}2)\Gamma(\frac{1-s}2)Q_{m}(s)}
$$

\medskip

\textit{Proof}

By Mellin transform of the relation of Lemma 1, we get 
$$
\int_0^{+\infty} t^{ s-1}\frac 1{m!}\Psi_m(t)dt
=\Big(\int_0^{+\infty}2\sqrt{2}e^{-\pi/u^2}\frac1u \ u^{s-1}du\Big)
\Big(\int_0^{+\infty} \frac{1}{(2m)!}\Phi_{2m}(x)x^{s-1}dx\Big)
$$

that is
 $$
\int_0^{+\infty} t^{ s-1}\frac 1{m!}\Psi_m(t)dt
=\sqrt{2}\pi^{\frac{s-1}2}\Gamma(\frac{1-s}2)
\int_0^{+\infty}  \frac{1}{(2m)!}\Phi_{2m}(x)x^{s-1}dx
$$
$\square$

\bigskip

\textbf{Remark}

Using $$\frac1t\Psi_m(\frac1t)=(-1)^m\Psi_m(t)$$ 

we get with the change of variable $t\mapsto \frac1t$ 
$$
\int_0^{+\infty} t^{ s-1}\frac 1{m!}\Psi_m(t)dt
=
(-1)^m\int_0^{+\infty} t^{ -s}\frac 1{m!}\Psi_m(t)dt
$$

 for $0<Re(s)<1$.
 
 By the preceding lemma this gives
$$Q_m(1-s)=(-1)^mQ_m(s)$$

As a consequence of this relation we see that for  $s=\frac12+it$   the polynomials $t\mapsto Q_{2m}(\frac12+it)$ are in $\mathbb{R}[X]$.

\subsection{Expansion of Mellin transforms in terms of the polynomials $Q_m$}

If  we have for a function $f$ holomorphic in $S$ an expansion 
$$
f(t)=\sum_{m\geq 0}a_{2m}\frac {\Psi_{2m}(t)
}{({2m})!}$$
 
and if we can evaluate the Mellin transform of $f$ for  $0<Re(s)<1$ by  integration
of the terms of the series :
$$
\int_0^{+\infty}\Big(\sum_{m\geq 0}a_{2m}\frac {\Psi_{2m}(t)}{({2m})!}\Big) t^{s-1} dt =\sum_{m\geq 0}\frac {a_{2m}}{({2m})!}\int_0^{+\infty} t^{s-1}\Psi_{2m}(t) dt 
$$

then we get for  $0<Re(s)<1$
$$
{\int_0^{+\infty} f(t)t^{s-1}dt=\frac1{\sqrt{2\pi}}\Gamma(\frac{s}2)\Gamma(\frac{1-s}2)\sum_{m\geq 0}\frac{a_{2m}}{(2m)!}
Q_{2m}(s)}$$

A simple condition to justify this calculation is 
$$
 \sum_{m\geq 0} \frac{\vert a_{2m}\vert}{(2m)!}<+\infty
 $$
 
Since in this case we have for $0<Re(s)<1$
\begin{eqnarray*}
\int_0^{+\infty}\sum_{m\geq 0}\vert t^{s-1}a_{2m}\frac {\Psi_{2m}(t)}{({2m})!}\vert dt 
\leq\sum_{m\geq 0}\frac {\vert a_{2m}\vert}{({2m})!}\int_0^{+\infty} t^{Re(s)-1}\frac {\sqrt{2}}{\sqrt{1+t^2}}dt<+\infty\end{eqnarray*}

\bigskip

\textbf{Example}

We have by relation (2) 

$$
 \frac1{1+t}=\frac1{{2}}\sum_{m\geq 0} a_{2m}\frac {\Psi_{2m}(t)}{(2m)!}
\text{ with }\ 
a_{2m}=\frac{(4m)!}{2^{4m}(2m+1)!}$$

By Stirling formula we have $ \frac{ a_{2m}}{(2m)!}=O(m^{-\frac32})$ thus $\sum_{m\geq 0} \frac{\vert a_{2m}\vert}{(2m)!}<+\infty$.

 And for $0<Re(s)<1$ we get   
$$
\frac{\pi}{\sin(\pi s)}=
\frac1{2\sqrt{2\pi}} \Gamma(\frac{s}2)\Gamma(\frac{1-s}2)
\sum_{m\geq 0} \frac{(4m)!}{(2m)!(2m+1)!}{2^{-4m}}Q_{2m}(s)
$$

\subsection{A conjecture for an expansion of Zeta in the critical strip}
For $0<Re(s)<1$ it is known (cf. [3]) that
    the Mellin transform of the function $t\mapsto G(t)-1-\frac1{t}$ is 
$$
\int_0^{+\infty} t^{s-1}(G(t)-1-\frac1{t})dt
=\Gamma(\frac s2)\pi^{-s/2}\zeta(s)
$$

We have seen in 2.1 that
$${
G(t)-1-\frac1t= \sum_{m\geq 0}{\alpha_{2m}}\Psi_{2m}(t)
\ \text{ with }\ 
\alpha_{2m}=S_{4m}-\frac{{2^{-4m+1}}(4m)!}{(2m)!(2m)!}}
$$

If we proceed by integration of the terms of the preceding series  we get
$$
\Gamma(\frac s2)\pi^{-s/2}\zeta(s)=\frac1{\sqrt{2 \pi}}\Gamma(\frac{s}2) \Gamma(\frac{1-s}2) \sum_{m\geq 0} { \alpha_{2m}}Q_{2m}(s)
$$

Unfortunately it seems that in this case $\sum_{m\geq 0} {\vert \alpha_{2m}\vert}=+\infty$ and the justification of the preceding section does not work.

\bigskip

\textbf{Conjecture}

For  $0<Re(s)<1$ the evaluation of  the Mellin transform of 
$G(t)-1-\frac1{t}$ by integration of the terms of the preceding series is valid and we get
$${
\zeta(s)=\frac1{\sqrt{2\pi}}\pi^{\frac{s}2}\Gamma(\frac{1-s}2) \sum_{m\geq 0} { \alpha_{2m}}\ Q_{2m}(s)}
$$

$$\text{ with }
\alpha_{2m}=S_{4m}-\frac{{2^{-4m+1}}(4m)!}{(2m)!(2m)!}$$

\bigskip

As we have seen the polynomials  
  $$Q_{2m}(s)=\ _2F_1(-2m,s/2;1/2;2)
$$
 are related to Mellin transforms of the Hermite functions   $\Phi_{4m} $ and they have their roots on the line $Re(s)=1/2$.

\bigskip

\textbf{Remarks}

 \textbf{1)} For the Riemann-Hardy function (cf. [3]) defined for  $t\in\mathbb{R}$ by
$$Z(t)=\pi^{\frac{-it}2}\frac{\Gamma(\frac14+i\frac t2)}{\vert \Gamma(\frac14+i\frac t2)\vert}\zeta(\frac12+it)
$$

the preceding conjecture gives
$$Z(t)=\frac1{\sqrt{2\pi}} \sum_{m\geq 0}  \alpha_{2m}f_{2m}(t)
$$

where the functions 
$$
f_{2m}(t)=\pi^{\frac{1}4}\vert \Gamma(\frac14+i\frac t2)\vert Q_{2m}(\frac12+it)
$$

 are  orthogonal  in $L^2(\mathbb{R})$.

\pagebreak

\textbf{2)} 
Other expressions of $\zeta(s)$ in the critical strip are obtained by the use of the MŸntz formula (cf. [7]): for  a  continuously differentiable function $F$ on $[0,+\infty[$ such that $F$ and $F'$  are $O(x^{-a})$, $(a>1)$ when $x\rightarrow \infty$, we have for $0<Re(s)<1$ 
$$
\zeta(s)\mathcal{M}F(s)=\mathcal{M}\Big(\sum_{n\geq 1} F(nt)-\frac1t\int_{0}^{+\infty}F(x)dx\Big)(s)
$$

We now show  that,  with  our preceding method, we can obtain  a simple  expansion of Zeta in the critical strip by applying this formula to the function  $F(x)=e^{-2\pi x}$.

For 
 $0<Re(s)<1$ we have 
\begin{equation*}
{(2\pi)^{-s}}{\Gamma(s)}\zeta(s)=\mathcal{M}(\sum_{n\geq 1}e^{- 2\pi nt}-\frac1{2\pi t})
=\int_0^{+\infty}f(t) t^{s-1}dt
\end{equation*}

where $$f(t)=\frac  {1}{e^{ 2\pi t}-1}-\frac1{ 2\pi t}$$

It is possible to get an expansion of  $f(t)$ using  the  Laguerre functions 
$$
\varphi_m(x)= e^{- x}L_m(2 x)
$$

defined by the generating function 
$\frac1{1-u}e^{- x\frac {1+u}{1-u}}=\sum_{m\geq 0}e^{- x}L_m(2x)u^m.
$
\medskip

These functions are orthogonal in $L^2(]0,+\infty[)$ and if $J_0$ is the Bessel function of order 0 then (cf.[6])
$$
\int_0^{+\infty} J_0(2 \sqrt{\xi x})\varphi_m(x)dx=(-1)^m\varphi_m(\xi)
$$

For $t>0$ we set $$
\psi_m(t)= (\frac{t-1}{t+1})^m\frac2{{1+t}}
$$

Using the generating function of the $\varphi_m$ 
 we get 
\begin{equation}
 e^{-2\pi x t}= \sum_{m\geq 0}\varphi_{m}(2\pi x) \psi_m(t)
\end{equation}

Summing (3)
  for $x=n\geq 1$ we have formally for $Re(t)>0$
$$
\frac  {1}{e^{ 2\pi t}-1}=\sum_{m\geq 0}s_{m}\psi_m(t)
\text{ with } s_m= \sum_{n\geq 1} \varphi_m(2\pi n)
$$

Since 
$\frac 1{2\pi t} =\frac 1{2\pi }\sum_{m\geq 0}(-1)^m{\psi_m(t)}
$
 we get  for $Re(t)>0$
\begin{equation}
f(t)=\sum_{m\geq 0}\big(s_{m}-(-1)^m \frac 1{2\pi}\big)\psi_m(t)
\end{equation}

\medskip
The functions $\psi_m$ and  $\varphi_m$ are related  by
 the multiplicative convolution
\begin{equation}
\psi_m=2(-1)^m (e^{-\frac{1}x}\frac{1}x) * \varphi_m
\end{equation}

The Mellin transform of $\varphi_m$ is (cf. [4]) for $Re(s)>0$
$$
\int_0^{+\infty} \varphi_m(x)x^{s-1}dx
=\Gamma(s)q_m(s)$$

where $q_m$ is the polynomial 
$q_m(s)=\  _2F_1(-m,s;1;2)$.

By the orthogonality relation of the $\varphi_m$ we deduce that the polynomials $t\mapsto q_m(\frac12+it)$ are orthogonal  with respect to the Borel measure $\vert \Gamma(\frac12+i\frac t2)\vert^2 dt $ on $]0,+\infty[$. Thus $q_m$ has  his roots on the line $Re(s)=\frac12$.

By  (5) for $0<Re(s)<1$ we have the Mellin transform of $\psi_m$
$$
\mathcal{M}( \psi_m)(s)=2\Gamma(s)\Gamma(1-s)(-1)^mq_m(s) 
$$

By Mellin transform of (4) we get formally
$$
\zeta(s)=2{(2\pi)^{s}} \Gamma(1-s)
\sum_{m\geq 0}\Big((-1)^ms_m-\frac 1{2\pi}\Big)q_m(s)
$$

\bigskip

A more simple expansion can be obtained using the Mellin transform of 
the function $g(t)=\frac1t f(\frac1t)$.

Using the  Poisson formula we have
$$
g(t)=\frac1t\sum_{n\geq 1} e^{-2\pi\frac{ n}t}- \frac1{2\pi}
=\frac1{\pi}\sum_{n\geq 1} \frac1{1+n^2t^2}-\frac1{2t}$$

By MŸntz formula we get
$$
\mathcal{M}(g)(s)=\mathcal{M}\Big(\frac1{\pi}\sum_{n\geq 1} \frac1{1+n^2t^2}-\frac1{2t}\Big)(s)=\zeta(s)\frac1{2\pi}{\Gamma(\frac s2)} {\Gamma(1-\frac s2)} 
$$

We now apply our preceding method, to the function  $F(x)=\frac1{\pi}\frac1{1+x^2}$.

We verify  immediately  that  
\begin{equation}
\frac1{\pi}\frac1{1+x^2t^2}=\frac1{2\pi}\sum_{m\geq 0} (-1)^m\psi_m(x^2)\psi_m(t^2) 
\end{equation}

We have for $t>0$
$$
\frac1{\pi}\sum_{n\geq 1} \frac1{1+n^2t^2}=\frac1{2\pi}\sum_{m\geq 0} (-1)^m\sigma_m\psi_m(t^2)
\ \text{ where } \ \sigma_m=\sum_{n\geq 1}\psi_m(n^2)
$$

With $u=\frac{t^2-1}{t^2+1}$ we have $\frac1t=\frac2{1+t^2}(1-u^2)^{-1/2}$, thus we get
 $$\frac1t=\sum_{m\geq 0} c_m\psi_m(t)  
\text{ with }c_{2n}=\frac{(2n)!}{2^{2n}(n!)^2} \ \text{ and }c_{2n+1}=0$$

Finally  we have the expansion
\begin{equation}
g(t)=\frac1{2\pi}\sum_{m\geq 0}(-1)^m(\sigma_m-\pi c_m
)\psi_m(t^2)
\end{equation}

By Mellin transform of (7) we get formally
$$
 \zeta(s)=\sum_{m\geq 0} \Big(\sigma_m-\pi c_m
\Big)q_m(\frac s2) 
$$

Note that, unlike the preceding expansions related to Hermite and Laguerre functions, in this expansion the sequence 
$$m\mapsto \sigma_m-\pi c_m=\sum_{n\geq 1}(\frac{n^2-1}{n^2+1})^m\frac2{{1+n^2}}-\pi c_m$$  has a very regular oscillation with amplitude near $\frac{\sqrt{\pi}}{\sqrt{2}}\frac1{\sqrt{ m}}$, but the polynomials $s\mapsto  q_m(\frac s2) $ have their roots on the line $Re(s)=1$.

\bigskip
\pagebreak

\textbf{References}

\medskip

[1] D.Bump, K.K.Choi, P.Kurlberg, J.Vaaler. A local Riemann hypothesis.  Math. Zeitschrift 233. (2000).

[2] D.Bump , E.K.-S.Ng. On the Riemann Zeta Function. Math. Zeitschrift 192. (1986).

[3] H.M.Edwards. Riemann Zeta function. Dover. (1974).

[4] I.S.Gradshteyn, I.M.Ryzhik. Tables of Integrals, Series and Products. Academic Press, Inc. (1994).

[5] E.Hille. A Class of Reciprocal Functions. Annals of Maths. Second Series. Vol. 27  N¡4 (1926)

[6] N.N.Lebedev. Revised and translated by R.A.Silverman. Special functions and their applications. Dover (1972).

[7] E.C. Titchmarsh. D.R. Heat-Brown. The theory of the Riemann Zeta  function. Clarendon 1986.

\bigskip

\textbf{Acknowledgments.} 

My warmest thanks go to F. Rouvire, M. Miniconi and J.F. Burnol for their interest and helpful comments.

\bigskip

\section{Appendix. Another expression for the constants $S_{4m}$}

Using   Poisson summation formula we deduce that for $\varphi\in\mathcal{S}(\mathbb{R})$ 

$$
\sum_{(k,l)\in \mathbb{Z}^{2}}\int_{\mathbb{R}}e^{-2i\pi xk}
e^{-\pi(x-l)^2}\varphi(x)dx=\sum_{l\in  \mathbb{Z}}e^{-\pi l^2} \sum_{k\in \mathbb{Z}}\varphi (k)=S_0\sum_{k\in \mathbb{Z}}\varphi (k)
$$

Taking  $u\in \mathbb{C}$ and  $\varphi(x)= e^{-2\pi xu}e^{-\pi x^2}$, we have
$$
\int_{\mathbb{R}}e^{-2i\pi xk}
e^{-\pi(x-l)^2}\varphi(x)dx=\frac1{\sqrt{2}}(-1)^{kl}e^{-\pi\frac12(k^2+l^2)}
 e^{i\pi u(k+il)}e^{\pi  u^2/2}
 $$
 
This gives for  $u\in\mathbb{C}$ the relation
\begin{equation}
\sum_{n\in \mathbb{Z}}e^{-\pi n^2+2\pi nu-\pi u^2/2}=\frac{1
}{S_0\sqrt{2}}\sum_{(k,l)\in  \mathbb{Z}^2}(-1)^{kl}e^{-\pi\frac12(k^2+l^2)}e^{ i\pi u(k+il)}\end{equation}

Let us now define for every integer $m\geq 0$ 
$${
T_m=\sum_{(k,l)\in  \mathbb{Z}^2}(-1)^{kl}e^{-\pi\frac12(k^2+l^2)}(k+il)^m}$$

We have clearly $T_{2m+1} =0$ since
$$
(-1)^{-k(-l)}e^{-\pi\frac12((-k)^2+(-l)^2)}(-k-il)^{2m+1}=-[(-1)^{kl}e^{-\pi\frac12(k^2+l^2)}(k+il)^{2m+1}]
$$
 
thus $T_{4m+1} =T_{4m+3} =0$ and also $T_{4m+2} =0$ because
$$
(-1)^{(-k)l}e^{-\pi\frac12(l^2+(-k)^2)}(-k+il)^{4m+2}=-[(-1)^{kl}e^{-\pi\frac12(k^2+l^2)}(l+ik)^{4m+2}]
$$

Thus only the constants  $T_{4m} $ are non zero, and by  derivation with respect to $u$ of the holomorphic function defined by the right side of (6) we have
$$
\sum_{n\in \mathbb{Z}}e^{-\pi n^2+2\pi nu-\pi u^2/2}=\frac{1}{S_0\sqrt{2}}\sum_{m\geq 0}\pi^{4m}T_{4m}\frac{u^{4m}}{(4m)!}
$$

Now using the generating function of Hermite polynomials we have 

$$
e^{-\pi n^2+2\pi nu-\pi{u^2}/2}=\sum_{m\geq 0} (\frac{\pi}2 )^{\frac m2}\Phi_m (n)\frac{u^m}{m!}
$$

By summation with $n\in \mathbb{Z}$ of this relation  we deduce that for $\vert u \vert <1$

$$
\sum_{n\in \mathbb{Z}}e^{-\pi n^2+2\pi nu-\pi u^2/2}=\sum_{m\geq 0} (\frac{\pi}2 )^{\frac m2}\Big(\sum_{n\in \mathbb{Z}} \Phi_{m}(n)\Big) \frac{u^m}{m!}
$$

(the interchange of $\sum_{n\in \mathbb{Z}}$ and $\sum_{m\geq 0}$ is easily justified using Lemma 0).

Thus we have for $\vert u \vert <1$
$$
\sum_{m\geq 0} (\frac{\pi}2 )^{\frac m2}\Big(\sum_{n\in \mathbb{Z}} \Phi_{m}(n)\Big) \frac{u^m}{m!}=
\frac{1}{S_0\sqrt{2}}\sum_{m\geq 0}\pi^{4m}T_{4m}\frac{u^{4m}}{(4m)!}
$$

and by  identification we get

$$ 
\sum_{n\in \mathbb{Z}} \Phi_{4m+2}(n)=0\ \text{ and }\ 
\sum_{n\in \mathbb{Z}} \Phi_{4m}(n)=
\frac
{(2\pi)^{2m}}{S_0\sqrt{2}}T_{4m} 
$$

\end{document}